\documentclass[runningheads]{llncs}

\usepackage{graphicx}
\usepackage{graphics}
\usepackage{amssymb, latexsym, amsmath, epsfig}
\usepackage{color}
\usepackage{epstopdf}
\usepackage{soul}
\usepackage[normalem]{ulem}
\usepackage{float}
\usepackage{epsfig}
\usepackage{subfigure}

\graphicspath{{figures/}}

\newcommand{\bfs}[1]{{\boldsymbol #1}}

\begin{document}
\title{Residual minimization for isogeometric analysis in reduced and mixed forms\thanks{This work was partially supported by the European Union's Horizon 2020 Research and Innovation Program of the Marie Sk{\l}odowska-Curie grant agreement No. 777778, the Mega-grant of the Russian Federation Government (N 14.Y26.31.0013), the Institute for Geoscience Research (TIGeR), and the Curtin Institute for Computation.}}
%
%

\author{Victor M. Calo\inst{1,2} \and
Quanling Deng\inst{1} \and
Sergio Rojas\inst{1} \and 
Albert Romkes\inst{3}
}
\authorrunning{V. Calo et al.}

\institute{Department of Applied Geology, Curtin University, Kent Street, Bentley, Perth, WA 6102, Australia \\\email{Victor.Calo@curtin.edu.au,Quanling.Deng@curtin.edu.au, srojash@gmail.com} \and
Mineral Resources, Commonwealth Scientific and Industrial Research Organisation (CSIRO), Kensington, Perth, WA 6152, Australia \and
Department  of  Mechanical  Engineering,  South  Dakota  Schoolol  Mines  \&  Technology,  501  E.  St.Joseph  Street,  Rapid  City,  SD  57701,  USA \email{Albert.Romkes@sdsmt.edu}}

\maketitle              
\begin{abstract}
Most variational forms of isogeometric analysis use highly-continuous basis functions for both trial and test spaces. For a partial differential equation with a smooth solution, isogeometric analysis with highly-continuous basis functions for trial space results in excellent discrete approximations of the solution. However, we observe that high continuity for test spaces is not necessary. In this work, we present a framework which uses highly-continuous B-splines for the trial spaces and basis functions with minimal regularity and possibly lower order polynomials for the test spaces. To realize this goal, we adopt the residual minimization methodology. We pose the problem in a mixed formulation, which results in a system governing both the solution and a Riesz representation of the residual. We present various variational formulations which are variationally-stable and verify their equivalence numerically via numerical tests. 

\end{abstract}

\section{Introduction} \label{sec:intr} 

Isogeometric analysis, introduced in 2005~\cite{hughes2005isogeometric,cottrell2009isogeometric}, is a widely-used numerical method for solving partial differential equations (PDEs). This analysis unifies the finite element methods with computer-aided design tools. Within the framework of the classic Galerkin finite element methods, isogeometric analysis uses as basis functions the functions that describe the geometry in the computer-aided design (CAD) and engineering (CAE) technologies, namely, B-splines or non-uniform rational B-splines (NURBS)  and their generalizations. These CAD/CAE basis functions may possess higher continuity. This smoothness improves the numerical approximations of PDEs which have highly regular solutions. Bazilevs \textit{et al.} established the approximation, stability, and error estimates in \cite{bazilevs2006isogeometric}.  Cottrell \textit{et al.}~in \cite{cottrell2006isogeometric} first used isogeometric analysis to study the structural vibrations and wave propagation problems. Spectral analysis shows that isogeometric elements significantly improve the accuracy of the spectral approximation when compared with classical finite elements. In~\cite{hughes2014finite}, the authors explored the additional advantages of isogeometric analysis on the spectral approximation over finite elements. Moreover, the work~\cite{deng2018dispersionm,puzyrev2017dispersion,calo2017dispersion,deng2018dispersion} minimized the spectral approximation errors for isogeometric elements. The minimized spectral errors possess superconvergence (two extra orders) as the mesh is refined.

Galerkin isogeometric analysis uses highly-continuous B-splines or NURBS for both trial and test spaces. For PDEs with smooth solutions (high regularities), isogeometric elements based on smooth functions render solutions which have a better physical interpretation. For example, for a PDE with a solution in $H^2$, both the exact solution field and exact flux field are expected to be globally continuous. Isogeometric analysis with $C^1$ B-spline basis functions produces a globally continuous flux field while the classical finite element method does not. Therefore, for smooth problems, it makes sense to use highly-continuous B-splines for the trial spaces. However, it is not necessary to use highly-continuous B-splines for the test spaces. In the view of the work~\cite{collier2012cost} which explores the cost of continuity, we expect an extra cost for solving the resulting system when we apply unnecessary high continuities for the basis functions in the test spaces. This extra cost per degree of freedom for highly-continuous discretizations was verified for direct \cite{collier2012cost,pardo2012survey,collier2014computational} and iterative solvers \cite{collier2013cost}. Additionally, the work  \cite{garcia2017value,garcia2018refined} shows that a reduction in the continuity of the trial and test spaces may lead to faster and more accurate solutions. Thus, we seek to develop an isogeometric analysis which uses minimal regularity for the test spaces. In this work, we establish a framework which uses highly-continuous B-splines for the trial spaces and basis functions with minimal regularity for the test spaces.

To realize this goal, we adopt the residual minimization methodology sharing with the Discontinous Petrov-Galerkin (DPG) method, the idea of a Riesz representation for the residual with respect to a stabilizing norm (c.f., \cite{demkowicz2010class,demkowicz2011class,demkowicz2012class,zitelli2011class}). The main idea of DPG is the use of discontinuous or broken basis functions for the test spaces within the Petrov-Galerkin framework. The method computes the optimal test function by using a trial-to-test operator.  The goal is to automatically stabilize the discrete formulation \cite{niemi2013automatically,calo2014analysis,niemi2011discontinuous,niemi2013discontinuous} without parameter tuning. 
The DPG method can be interpreted as a minimum residual method in which the residual is measured in the dual test norm~\cite{demkowicz2017discontinuous}. By introducing an auxiliary unknown $\phi$ representing the Riesz representation of residual, the method is cast into a mixed problem. Effectively, we extend to highly continuous isogeometric discretizations the method described in \cite{calo2018automatic} where standard $C^0$ finite element basis functions build the trial space while the test space uses broken polynomials.

The proposed method has a feature which distinguishes us from DPG. We use highly-continuous B-splines for space $V^h$ in which we seek for the solution, while DPG uses discontinuous basis functions. Consequently, DPG introduces additional unknowns that live on the mesh skeleton; see equation (10) or (39) in~\cite{demkowicz2017discontinuous}.  We do not introduce any additional unknowns. Instead, we use basis  functions ($C^{-1}$ functions in the discrete $L^2$ space or at most $C^0$ in the discrete $H^1$ space defined on the mesh partition) with minimal regularity in the sense that no inner products are introduced on the mesh skeleton (thus, the bilinear form only involves element integrations) while maintaining the inf-sup condition for the resulting variational formulation.

The rest of this paper is organized as follows. Section \ref{sec:ps} describes the problem under consideration and introduces various variational formulations following the DPG framework closely. Section 3 presents isogeometric formulations at the discrete level, while Section 4 shows numerical examples to demonstrate the performance of the formulations.  Concluding remarks are given in Section 5.

\section{Problem statement} \label{sec:ps}

Let $\Omega \subset \mathbb{R}^d, d=1,2,3,$ be an open bounded domain with Lipschitz boundary $\partial \Omega$. We consider the advection-diffusion-reaction equation: Find $u$ such that 
\begin{equation} \label{eq:pde}
\begin{aligned}
-\nabla \cdot (\kappa \nabla u - \bfs{\beta} u) + \gamma u & = f \quad \text{in} \quad \Omega, \\
u & = 0 \quad \text{on} \quad \partial \Omega, \\
\end{aligned}
\end{equation}
where  $\kappa$ is the diffusion coefficient, $\bfs{\beta}$ is the advective velocity vector, $\gamma \ge 0$ is the reaction coefficient, $u$ is the function to be found, and $f$ is a forcing function. This problem can be cast into an equivalent system of two first-order equations
\begin{equation} \label{eq:pdee}
\begin{aligned}
\bfs{q} - \kappa \nabla u + \bfs{\beta} u & = \bfs{0} \quad \text{in} \quad \Omega,  \\
-\nabla \cdot \bfs{q} + \gamma u & = f \quad \text{in} \quad \Omega, \\
u & = 0 \quad \text{on} \quad \partial \Omega, \\
\end{aligned}
\end{equation}
where $\bfs{q}$ is an auxiliary variable standing for the flux.

\subsection{General setting} \label{sec:setting}
We present the general setting for equations \eqref{eq:pde} and \eqref{eq:pdee} as they have distinguishing features. Let $V$ and $W$ denote the trial and test space, resepectively. For \eqref{eq:pdee} we let $(u, \bfs{q})$ be the solution pair in its trial space $V^u \times V^{\bfs{q}}$ and let its corresponding test space be $W^u \times W^{\bfs{q}}$. We specify these spaces in each of the variational formulations. Let $L^2(\Omega) = \{ v: \int_\Omega v^2 \ \text{d} x < \infty \}$ and denote the induced $L^2$ norm as $\| \cdot \|_{0,\Omega}$.  Let $C^k(\Omega), k=-1,0,1, \cdots, \infty,$ denote the space of functions having continuous (partial) derivatives up to the $k$-th order over the whole domain $\Omega$.  In particular, for $k=-1$, we mean that the function in $C^{-1}$ is discontinuous somewhere in $\Omega$. Let $H^1(\Omega) = \{ v: v \in L^2(\Omega), \int_\Omega (\nabla v)^2 \ \text{d} x < \infty \}$ and $H^1_0(\Omega)$ be the space of all functions in $H^1(\Omega)$ vanishing at the boundary $\partial \Omega$. The $H^1$ seminorm is defined as $|w|_{1,\Omega} = \| \nabla w \|_{0,\Omega}$ and $H^1$ norm is defined as $\|w\|^2_{1,\Omega} = \| \nabla w \|^2_{0,\Omega} + |w|^2_{1,\Omega}$.  Finally, $H(\text{div}, \Omega)$ consists of all square integrable vector-valued fields on $\Omega$ whose divergence is a function that is also square integrable. Similarly, the $H(\text{div}, \Omega)$ norm is defined as $\| \bfs{p} \|^2_{\text{div}, \Omega} = \| \bfs{p} \|^2_{0,\Omega} + \|\nabla \cdot \bfs{p} \|^2_{0,\Omega}$.

Let $\mathcal{T}_h$ be a partition of $\Omega$ into non-overlapping mesh elements. For simplicity and the purpose of using B-splines in multiple dimensions, we assume the tensor-product structure.  Let $K\in \mathcal{T}_h$ be a generic element and denote by $\partial K$ its boundary. Let $\bfs{n}$ be the outward unit normal vector.
Let $(\cdot, \cdot)_S$ denote the $L^2(S)$ the inner product where $S$ is a $d$ or $d-1$ dimensional domain ($S$ is typically $\Omega, K, \partial \Omega, \partial K$). 

At discrete level, we define the finite spaces, namely the finite subspaces of $V, W, V^u \times V^{\bfs{q}}$, and $W^u \times W^{\bfs{q}}$. Since isogeometric analysis adopts highly-continuous basis functions, we specify the corresponding finite spaces with both the polynomial order $p$ and the order $k$ of global continuity (that is, $C^k$). Let us denote by $V^h_{p,k}, W^h_{p,k}$ and $V^{u,h}_{p,k}, V^{\bfs{q},h}_{p,k}, W^{u,h}_{p,k}, W^{\bfs{q},h}_{p,k}$ the finite-dimensional  subspaces of $V, W$ for \eqref{eq:pde} and $V^u, V^{\bfs{q}}, W^u, W^{\bfs{q}}$ for \eqref{eq:pdee}, respectively. Consequently, $V^h_{p,k} \subset V, W^h_{p,k} \subset W, V^{u,h}_{p,k} \subset V^u, V^{\bfs{q},h}_{p,k} \subset V^{\bfs{q}}, W^{u,h}_{p,k} \subset W^u$, and $W^{\bfs{q},h}_{p,k} \subset W^{\bfs{q}}$. In this work, we construct the test spaces $W^h_{q,l} $ and $W^{u,h}_{q,l}$ using basis functions with lower order polynomials as well as lower regularity.

\subsection{Various variational formulations at continuous level}
In this section, we present seven variational formulations for both \eqref{eq:pde} and \eqref{eq:pdee} following closely the six formulations described within the DPG framework  in~\cite{demkowicz2017discontinuous}. In particular, we add one more formulation (to the six formulations in DPG framework) which resembles the isogeometric collocation method.

We start with the abstract variational formulations of \eqref{eq:pde} and \eqref{eq:pdee}. The variational weak formulations of \eqref{eq:pde} can be written as: Find $u \in V$ such that 
\begin{equation} \label{eq:pvf}
b(w, u) = \ell(w), \qquad   \forall \ w \in W,
\end{equation}
where $b(\cdot, \cdot)$ and $\ell(\cdot)$ are bilinear and linear forms, respectively. Similarly, the variational weak formulations of \eqref{eq:pdee} can be written: Find $(u, \bfs{q}) \in V^u \times V^{\bfs{q}}$ such that 
\begin{equation} \label{eq:mvf}
b( (w, \bfs{p})  ,(u, \bfs{q})) = \ell( (w, \bfs{p}) ), \qquad   \forall \ (w, \bfs{p}) \in W^u \times W^{\bfs{q}},
\end{equation}
where $b(\cdot, \cdot)$ and $\ell(\cdot)$ are bilinear and linear forms defined over two fields. Herein, we refer to $(u, \bfs{q})$ as the solution (trial) pair while to $(w, \bfs{p})$ as the weighting (test) function pair. Equations \eqref{eq:pde} and \eqref{eq:pvf} are the primal forms of the PDE and the variational formulation while equations \eqref{eq:pdee} and \eqref{eq:mvf} are their mixed forms. We keep this structure for the forms at discrete level in Section \ref{sec:m}. All these forms and their associated spaces are to be specified in each particular formulation as follows.

\begin{enumerate}

\item Primal trivial formulation: Let $V = H^1_0(\Omega) \cap C^1(\Omega), W = L^2(\Omega)$. Find $u \in V $ satisfying \eqref{eq:pvf}  with 
\begin{equation} \label{eq:pf1}
\begin{aligned}
b(w,u) & := b_1( w,u) = (w,-\nabla \cdot (\kappa \nabla u - \bfs{\beta} u) + \gamma u)_\Omega, \\
\ell(w) &:= \ell_1(w) = (w, f)_\Omega.
\end{aligned}
\end{equation}

\item Primal classical (FEM) formulation: Let $V = H^1_0(\Omega), W = H^1_0(\Omega)$. Find $u \in V $ satisfying \eqref{eq:pvf}  with 
\begin{equation} \label{eq:pf2}
\begin{aligned}
b(w,u) & := b_2(w,u) = (\nabla w,  \kappa \nabla u - \bfs{\beta} u)_\Omega + (w, \gamma u )_\Omega, \\
\ell(w) &:= \ell_2(w) = (w, f)_\Omega.
\end{aligned}
\end{equation}

\item  Mixed trivial formulation: Let $V^u = H^1_0(\Omega), V^{\bfs{q}} = H(\text{div}, \Omega), W^u = L^2(\Omega), $ $W^{\bfs{q}} = (L^2(\Omega))^d$. Find $(u, \bfs{q}) \in V^u \times V^{\bfs{q}}$ satisfying \eqref{eq:mvf} with 
\begin{equation} \label{eq:mf1}
\begin{aligned}
b((w, \bfs{p}), (u, \bfs{q}) ) & := b_3((w, \bfs{p}), (u, \bfs{q}) ) \\
& = (w, -\nabla \cdot \bfs{q} + \gamma u)_\Omega + (\bfs{p}, \bfs{q} - \kappa \nabla u + \bfs{\beta} u)_\Omega, \\
\ell( (w, \bfs{p}) ) &:= \ell_3( (w, \bfs{p}) ) = (w, f)_\Omega + (\bfs{0}, \bfs{q})_\Omega = (w, f)_\Omega.
\end{aligned}
\end{equation}

\item Mixed classical formulation: Let $V^u = L^2(\Omega), V^{\bfs{q}} = H(\text{div}, \Omega), W^u = L^2(\Omega), $ $ W^{\bfs{q}} = H(\text{div}, \Omega)$. Find $(u, \bfs{q}) \in V^u \times V^{\bfs{q}}$ satisfying \eqref{eq:mvf} with 
\begin{equation} \label{eq:mf2}
\begin{aligned}
b((w, \bfs{p}), (u, \bfs{q}) ) & := b_4((w, \bfs{p}), (u, \bfs{q}) ) \\ 
& = (w, -\nabla \cdot \bfs{q} + \gamma u)_\Omega + (\bfs{p}, \bfs{q} + \bfs{\beta} u)_\Omega + (\nabla \cdot (\kappa \bfs{p}), u )_\Omega, \\
\ell( (w, \bfs{p}) ) &:= \ell_4( (w, \bfs{p}) ) = (w, f)_\Omega + (\bfs{0}, \bfs{q})_\Omega = (w, f)_\Omega.
\end{aligned}
\end{equation}

\item Mixed classical formulation: Let $V^u = H^1_0(\Omega), V^{\bfs{q}} = (L^2(\Omega))^d, W^u = H^1_0(\Omega),$ $ W^{\bfs{q}} = (L^2(\Omega))^d$. Find $(u, \bfs{q}) \in V^u \times V^{\bfs{q}}$ satisfying \eqref{eq:mvf} with 
\begin{equation} \label{eq:mf3}
\begin{aligned}
b((w, \bfs{p}), (u, \bfs{q}) ) & := b_5((w, \bfs{p}), (u, \bfs{q}) ) \\
& = (\nabla w, \bfs{q})_\Omega + (w, \gamma u )_\Omega  + (\bfs{p}, \bfs{q} - \kappa \nabla u + \bfs{\beta} u)_\Omega, \\
\ell( (w, \bfs{p}) ) &:= \ell_5( (w, \bfs{p}) ) = (w, f)_\Omega + (\bfs{0}, \bfs{q})_\Omega = (w, f)_\Omega.
\end{aligned}
\end{equation}

\item Mixed ultraweak formulation: Let $V^u = L^2(\Omega), V^{\bfs{q}} = (L^2(\Omega))^d, W^u = H^1_0(\Omega), W^{\bfs{q}} = H(\text{div}, \Omega)$. Find $(u, \bfs{q}) \in V^u \times V^{\bfs{q}}$ satisfying \eqref{eq:mvf} with 
\begin{equation} \label{eq:mf4}
\begin{aligned}
b((w, \bfs{p}), (u, \bfs{q}) ) & := b_6((w, \bfs{p}), (u, \bfs{q}) ) \\
& = (\nabla w, \bfs{q})_\Omega + (w, \gamma u)_\Omega + (\bfs{p}, \bfs{q} + \bfs{\beta} u)_\Omega + (\nabla \cdot (\kappa \bfs{p}), u )_\Omega, \\
\ell( (w, \bfs{p}) ) &:= \ell_6( (w, \bfs{p}) ) = (w, f)_\Omega + (\bfs{0}, \bfs{q})_\Omega = (w, f)_\Omega.
\end{aligned}
\end{equation}

\item Reduced flux formulation: Let $V = H(\text{div}, \Omega), W = H(\text{div}, \Omega)$. Find $\bfs{q} \in V $ satisfying \eqref{eq:pvf}  with 
\begin{equation} \label{eq:rff}
\begin{aligned}
b(\bfs{p}, \bfs{q}) & := b_7(\bfs{p}, \bfs{q}) = (\bfs{p}, \bfs{q})_\Omega + (\nabla \cdot (\kappa \bfs{p}), \gamma^{-1} \nabla \cdot \bfs{q})_\Omega + (\bfs{p}, \bfs{\beta}\gamma^{-1} \nabla \cdot \bfs{q})_\Omega, \\
\ell(\bfs{p}) &:= \ell_7(\bfs{p}) = -(\nabla \cdot (\kappa \bfs{p}), \gamma^{-1}f)_\Omega - (\bfs{p}, \bfs{\beta}\gamma^{-1} f)_\Omega.
\end{aligned}
\end{equation}

\end{enumerate}
Herein, the primal trivial formulation reduces to the isogeometric collocation method when applying constant test functions.


\section{Various isogeometric formulations} \label{sec:m}
In this section, we present the isogeometric formulations at the discrete level for both \eqref{eq:pvf} and \eqref{eq:mvf}. We first specify the basis functions for all the finite element spaces associated with the mesh configuration $\mathcal{T}_h$. For this purpose, we use the Cox-de Boor recursion formula~\cite{de1978practical,piegl2012nurbs} on each dimension and then take tensor-product to obtain the necessary basis functions for multiple dimensions. The Cox-de Boor recursion formula generates $C^k$ and $p$-th order B-spline basis functions, where $p=0,1,2,\cdots,$ and $k=-1, 0,1,\cdots, p-1$. $C^{-1}$ basis functions generate a finite-dimensional subspace of the $L^2(\Omega)$ and $(L^2(\Omega))^d$, while $C^k, k=0,1,\cdots,p-1$ basis functions generate a finite-dimensional subspace of the $H^1(\Omega)$ and $H(\text{div}, \Omega)$.  These finite-dimensional subspaces consist of piecewise polynomials of order $p=0,1,2,\cdots$ and of continuity order $k=-1,0,1,\cdots,p-1$.

The definition of the B-spline basis functions in one dimension is as follows. 
Let $X = \{x_0, x_1, \cdots, x_m \}$ be a knot vector with knots $x_j$, that is, a nondecreasing sequence of real numbers which are called knots.  The $j$-th B-spline basis function of degree $p$, denoted as $\theta^j_p(x)$, is defined as~\cite{de1978practical,piegl2012nurbs}
\begin{equation} \label{eq:B-spline}
\begin{aligned}
\theta^j_0(x) & = 
\begin{cases}
1, \quad \text{if} \ x_j \le x < x_{j+1} \\
0, \quad \text{otherwise} \\
\end{cases} \\ 
\theta^j_p(x) & = \frac{x - x_j}{x_{j+p} - x_j} \theta^j_{p-1}(x) + \frac{x_{j+p+1} - x}{x_{j+p+1} - x_{j+1}} \theta^{j+1}_{p-1}(x).
\end{aligned}
\end{equation}

We then construct the finite-dimensional subspaces using the B-splines on uniform tensor-product meshes with non-repeating and repeating knots. For a $p$-th order B-spline, a repetition of $k =0,1,\cdots,p$ times of an internal node results in a function of $C^{p-1-k}$ continuity; see~\cite{piegl2012nurbs,cottrell2009isogeometric}. These B-spline basis functions characterize the finite-dimensional subspaces. For example,
\begin{equation} \label{eq:bsV}
\begin{aligned}
V^h_{p,k} & = 
\begin{cases}
S^p_k = \text{span} \{ \theta_j^p(x) \}_{j=1}^{N_x}, & \text{in 1D}\\
S^{p, p}_{k, k} = \text{span} \{ \theta_i^p(x) \theta_j^p(y) \}_{i, j=1}^{N_x, N_y}, & \text{in 2D}\\
S^{p, p,p}_{k, k,k} = \text{span} \{ \theta_i^p(x) \theta_j^p(y) \theta_l^p(z) \}_{i, j,l=1}^{N_x, N_y,N_z}, & \text{in 3D}\\
\end{cases} \\
V^{\bfs{q},h}_{p,k} & = 
\begin{cases}
S^{p-1}_{k-1}, & \text{in 1D}\\
S^{p-1, p-1}_{k-1, k-1} \times S^{p-1, p-1}_{k-1, k-1}, & \text{in 2D}\\
S^{p-1, p-1}_{k-1, k-1} \times S^{p-1, p-1}_{k-1, k-1} \times S^{p-1, p-1}_{k-1, k-1}, & \text{in 3D}\\
\end{cases} \\
\end{aligned}
\end{equation}
where $p$ and $k$ specify the approximation order and continuity order in each dimension (they can be different in general), respectively. $N_x, N_y, N_z$ is the total number of basis functions in each dimension. Note that the space $V^{\bfs{q},h}_{p,k}$ collapse to standard Raviart-Thomas mixed finite elements \cite{raviart1977mixed}; see \cite[Section 5]{evans2013isogeometric} or \cite[Section 3]{buffa2010isogeometric} for more details on the construction of these spaces using B-splines.

We now adopt the mixed formulation for \eqref{eq:pvf} at discrete level: Find $(u^h, \phi^h) \in V^h_{p,k} \times W^h_{q,l}$ such that 
\begin{equation} \label{eq:dpvf}
\begin{aligned}
g(w^h, \phi^h) + b(w^h, u^h) & = \ell(w^h), \qquad   \forall \ w^h \in W^h_{q,l}, \\
b(\phi^h, v^h) & = 0, \qquad  \qquad \forall \ v^h \in V^h_{p,k}, \\
\end{aligned}
\end{equation}
where the spaces are chosen such that $\text{dim}(W^h_{q,l}) \ge \text{dim}(V^h_{p,k})$. 
Herein, $\phi^h$ is the negative of the Riesz representation of the residual. The auxiliary bilinear form $g(\cdot, \cdot)$ is an inner product and it produces a Gramm matrix for the purpose of residual minimization. We define it generally as follows
\begin{equation} \label{eq:ra}
g(v, w)  = \sum_{K \in \mathcal{T}_h } \tau_0 (v, w)_K + \tau_1 h^{\iota_1} (\nabla v, \nabla w)_K + \tau_2 h^{\iota_2} (\Delta v, \Delta w)_K, \quad   \forall \ v, w \in W^h_{q,l}, \\
\end{equation}
where $\tau_i, \iota_j \in \mathbb{R}, i=0,1,2, j=1,2$ are free parameters.  The default setting is $\tau_0=1, \tau_1=1, \tau_2=0$, $\iota_1 = 2, \iota_2=0$. 
Once one specifies an inner product $g(\cdot, \cdot)$, then under inf-sup assumption on the discrete bilinear formulation in \eqref{eq:dpvf}, the approximate solution $u^h$ has a minimal error in the energy norm induced from $g(\cdot, \cdot)$; see, for example~\cite{demkowicz2017discontinuous,demkowicz2011class}.

Similarly, we present the mixed formulation for \eqref{eq:mvf} at discrete level: Find $(u^h, \bfs{q}^h) \in V^{u,h}_{p_1,k_1} \times V^{\bfs{q},h}_{p_2,k_2}$ and $(\phi^h, \bfs{\psi}^h) \in \times W^{w,h}_{q_1,l_1} \times W^{\bfs{q},h}_{q_2,l_2}$ such that 
\begin{equation} \label{eq:dmvf}
\begin{aligned}
g((w^h, \bfs{p}^h), (\phi^h, \bfs{\psi}^h) ) + b((w^h, \bfs{p}^h), (u^h, \bfs{q}^h)) & = \ell( (w^h, \bfs{p}^h) ), \forall \ (w^h, \bfs{p}^h) \in W^{w,h}_{q_1,l_1} \times W^{\bfs{q},h}_{q_2,l_2}, \\
b((\phi^h, \bfs{\psi}^h), (v^h, \bfs{r}^h)) & = 0, \qquad \quad \forall \ (v^h, \bfs{r}^h) \in V^{u,h}_{p_1,k_1} \times V^{\bfs{q},h}_{p_2,k_2}, \\
\end{aligned}
\end{equation}
where the spaces are chosen such that $\text{dim}(W^h_{q_1,l_1}) \ge \text{dim}(V^{u,h}_{p_1,k_1})$ and $\text{dim}(W^{\bfs{q},h}_{q_2,l_2}) \ge \text{dim}(V^{\bfs{q},h}_{p_2,k_2})$. Herein, the continuity orders corresponding to solution $u$ and flux $\bfs{q}$ can be different. Potentially, their polynomial orders can also be different. Similarly, $g((\cdot, \cdot), (\cdot, \cdot))$ is an inner product for the purpose of residual minimization. We define the Gramm product generally as follows, for $(v, \bfs{r}), (w, \bfs{p}) \in W^{w,h}_{q_1,l_1} \times W^{\bfs{q},h}_{q_2,l_2}$, 
\begin{equation} \label{eq:ma}
\begin{aligned}
g((v, \bfs{r}), (w, \bfs{p}) )  & = \sum_{K \in \mathcal{T}_h } \tau_3 (v, w)_K + \tau_4 h^{\iota_3} (\nabla v, \nabla w)_K \\
& \qquad \qquad + \tau_5 (\bfs{r}, \bfs{p})_K + \tau_6 h^{\iota_4} (\nabla \cdot \bfs{r}, \nabla \cdot \bfs{p})_K,
\end{aligned}
\end{equation} 
where $\tau_i, \iota_j \in \mathbb{R}, i=3,4,\cdots,6, j=3,4$ are free parameters.  The default setting is $\tau_3=1, \tau_4=1, \tau_5=1, \tau_6=1, \iota_3 = 2, \iota_4 = 2$. The formulation \eqref{eq:dmvf} is the same as \eqref{eq:dpvf} if we view $(u^h, \bfs{q}^h)$ as a solution pair. 

Within these formulations, once we specify all free parameters, the bilinear and linear forms, and space settings, we have a different method. We present the following discrete variational formulations

\begin{enumerate}

\item Let $b(\cdot, \cdot) = b_1(\cdot, \cdot), \ell(\cdot) = \ell_1(\cdot)$ defined in \eqref{eq:pf1}. Let $V^h_{p,k}$ consist of B-spline basis functions of continuity $C^k, k\ge1$ and $W^h_{q,l}$ consist of discontinuous basis functions.  The \textit{discrete primal trivial formulation} is:  Find $(u^h, \phi^h) \in V^h_{p,k} \times W^h_{q,l}$ satisfying \eqref{eq:dpvf}.

\item Let $b(\cdot, \cdot) = b_2(\cdot, \cdot), \ell(\cdot) = \ell_2(\cdot)$ defined in \eqref{eq:pf2}. Let $V^h_{p,k}$ and $W^h_{q,l}$  consist of B-spline basis functions of continuity at least $C^0$.  The \textit{discrete primal classical formulation} is:  Find $(u^h, \phi^h) \in V^h_{p,k} \times W^h_{q,l}$ satisfying \eqref{eq:dpvf}.

\item Let $b(\cdot, \cdot) = b_3(\cdot, \cdot), \ell(\cdot) = \ell_3(\cdot)$ defined in \eqref{eq:mf1}. Let the solution space consist of B-spline basis functions of continuity at least $C^0$ while the test space consist of discontinuous basis functions. The \textit{discrete mixed trivial formulation} is: Find $(u^h, \bfs{q}^h) \in V^{u,h}_{p_1,k_1} \times V^{\bfs{q},h}_{p_2,k_2}$ and $(\phi^h, \bfs{\psi}^h) \in \times W^{w,h}_{q_1,l_1} \times W^{\bfs{q},h}_{q_2,l_2}$ satisfying \eqref{eq:dmvf}.

\item Let $b(\cdot, \cdot) = b_4(\cdot, \cdot), \ell(\cdot) = \ell_4(\cdot)$ defined in \eqref{eq:mf2}. Let the solution space and test space for flux consist of B-spline basis functions of continuity at least $C^0$ while the test space for $u$ consist of discontinuous basis functions. The \textit{discrete mixed classical formulation I} is: Find $(u^h, \bfs{q}^h) \in V^{u,h}_{p_1,k_1} \times V^{\bfs{q},h}_{p_2,k_2}$ and $(\phi^h, \bfs{\psi}^h) \in \times W^{w,h}_{q_1,l_1} \times W^{\bfs{q},h}_{q_2,l_2}$ satisfying \eqref{eq:dmvf}.

\item Let $b(\cdot, \cdot) = b_5(\cdot, \cdot), \ell(\cdot) = \ell_5(\cdot)$ defined in \eqref{eq:mf3}. Let the solution space and test space for $u$ consist of B-spline basis functions of continuity at least $C^0$ while the test space for flux consist of discontinuous basis functions. The \textit{discrete mixed classical formulation II} is: Find $(u^h, \bfs{q}^h) \in V^{u,h}_{p_1,k_1} \times V^{\bfs{q},h}_{p_2,k_2}$ and $(\phi^h, \bfs{\psi}^h) \in \times W^{w,h}_{q_1,l_1} \times W^{\bfs{q},h}_{q_2,l_2}$ satisfying \eqref{eq:dmvf}.

\item Let $b(\cdot, \cdot) = b_6(\cdot, \cdot), \ell(\cdot) = \ell_6(\cdot)$ defined in \eqref{eq:mf4}. Let the solution and test spaces consist of B-spline basis functions of continuity at least $C^0$.  The \textit{discrete mixed ultraweak formulation II} is: Find $(u^h, \bfs{q}^h) \in V^{u,h}_{p_1,k_1} \times V^{\bfs{q},h}_{p_2,k_2}$ and $(\phi^h, \bfs{\psi}^h) \in \times W^{w,h}_{q_1,l_1} \times W^{\bfs{q},h}_{q_2,l_2}$ satisfying \eqref{eq:dmvf}.

\end{enumerate}

The primal classical formulation can be reduced to the standard finite element and isogeometric element methods.   If we set $k=0, l=0$ for $p\ge 1$, then we have $\phi^h = 0$ due to the orthogonal condition (second equation) in \eqref{eq:dpvf}. Thus, the method reduces to finite element method. Similarly, if we set $k=l=1, \cdots, p-1$ for $p\ge 2$, then we have $\phi^h = 0$ in \eqref{eq:dpvf} as well and the method reduces to isogeometric analysis. For all other discrete variational formulations, we may constrain the solution and test spaces in such a way that the variational forms we discuss render standard discretization techniques when the Riesz representation of the residual is identically zero. 

For the scenarios where we use different trial and test spaces, we obtain a non-zero discrete representation of the residual $\Phi$, which we use as an error estimator to guide the refinements of the meshes accordingly. The error estimators are defined in the sense of $G$ which is a result of the bilinear form $g(\cdot, \cdot)$. We refer to the DPG work \cite{demkowicz2010class,demkowicz2011class,demkowicz2012class,zitelli2011class} for more details in this direction.

These discrete variational formulations lead to a linear matrix problem
\begin{equation} \label{eq:mp}
\begin{bmatrix}
G & B \\
B^T & 0 \\
\end{bmatrix}
\begin{bmatrix}
\Phi \\
U \\
\end{bmatrix} 
= 
\begin{bmatrix}
L\\
\bfs{0} \\
\end{bmatrix},
\end{equation}
where $L$ is the corresponding forcing term arising from the linear form $\ell(\cdot)$, $\Phi, U$ are the solution pairs for the mixed formulations, $G$ represents the Gramm product in which we minimize the residual that arises from the bilinear form $g(\cdot, \cdot)$  and $B$ is the matrix arising from the bilinear form $b(\cdot, \cdot)$. We solve the first equation in \eqref{eq:mp} and substitute to the second equation in \eqref{eq:mp} to obtain
\begin{equation} \label{eq:rmp}
B^T G^{-1} ( L- B U ) = B^T G^{-1} L- B^T G^{-1} B U  = \bfs{0}.
\end{equation}
Herein, $B U -L$ is the residual and \eqref{eq:rmp} is a least-square type of problem.

For all the discrete variational formulations, we verify the following optimal error convergence rates numerically in the Section \ref{sec:num}:
\begin{equation} \label{eq:uh1err}
| u - u^h |_{1,\Omega} \le C h^p,
\end{equation}
where $C$ is a constant independent of $h$. For the approximate fluxes, we define the following errors in $L^2$ norm
\begin{equation}
\| \bfs{q} - \bfs{q}^h \|_{0,\Omega} = \| (\kappa \nabla u - \bfs{\beta} u) - \bfs{q}^h \|_{0,\Omega}.
\end{equation}
The optimal convergence rate is
\begin{equation} \label{eq:qhl2}
\| \bfs{q} - \bfs{q}^h \|_{0,\Omega} \le C h^{p+1},
\end{equation}
where $C$ is a constant independent of $h$.

\section{Numerical experiments} \label{sec:num}
The main result of these numerical tests is that the various formulations we discuss above are equivalent in the sense of resulting the same (optimal) error convergence rates. We focus on 2D and consider the problem \eqref{eq:pde} with  $\kappa = 1, \gamma=1, \bfs{\beta} = (1, 1)^T$ and a manufactured solution
\begin{equation}
u(x,y) = \sin(\pi x) \sin(\pi y) (2-x+3y).
\end{equation} 
$f$ is the corresponding forcing satisfying \eqref{eq:pde}. The true flux $\bfs{q}$ is calculated from \eqref{eq:pdee}. We apply all six variational formulations, namely, two primal and four mixed formulations, to solve this problem.

\begin{figure}[!ht]
\centering\includegraphics[width=6.0cm]{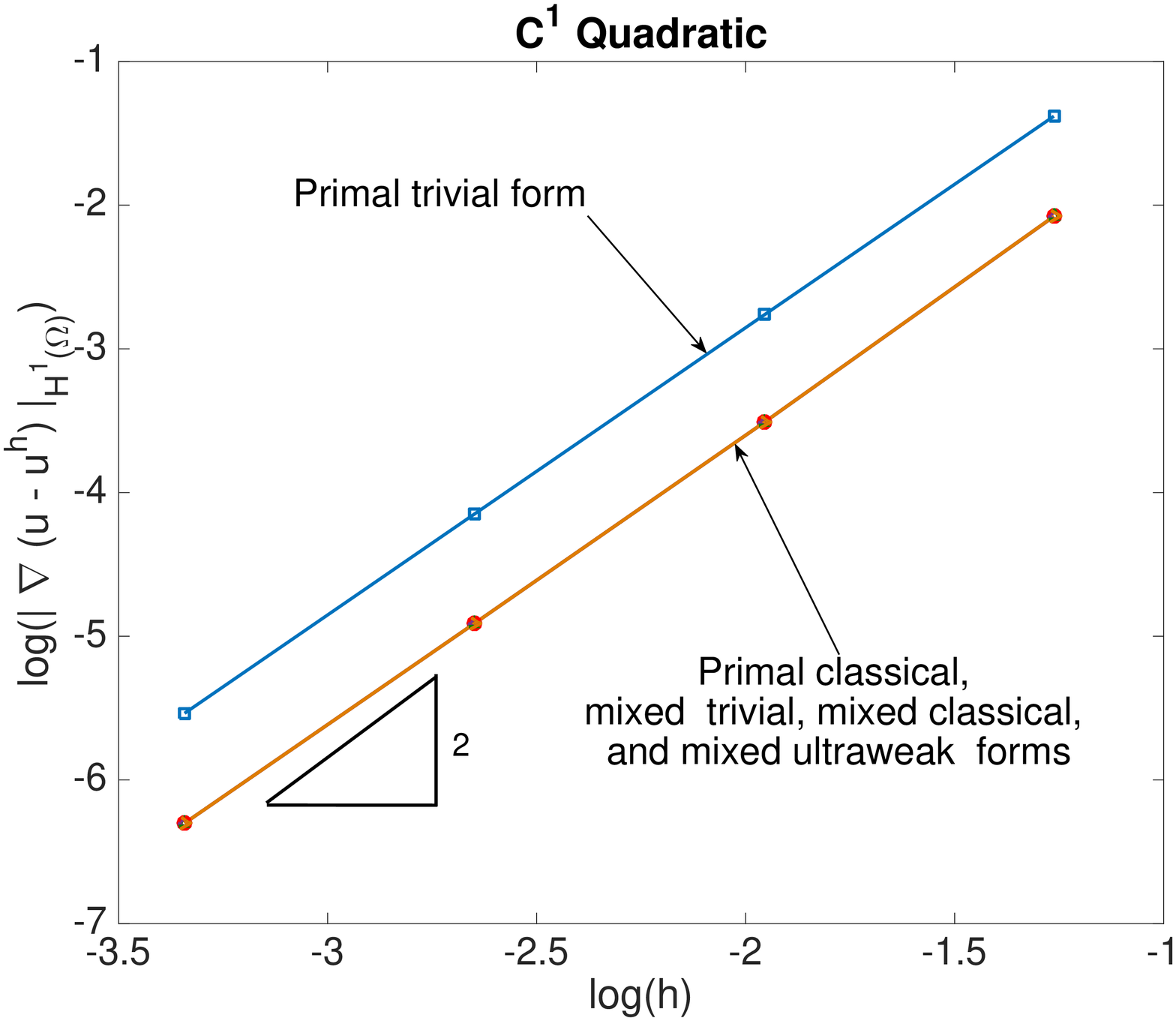} 
\centering\includegraphics[width=6.0cm]{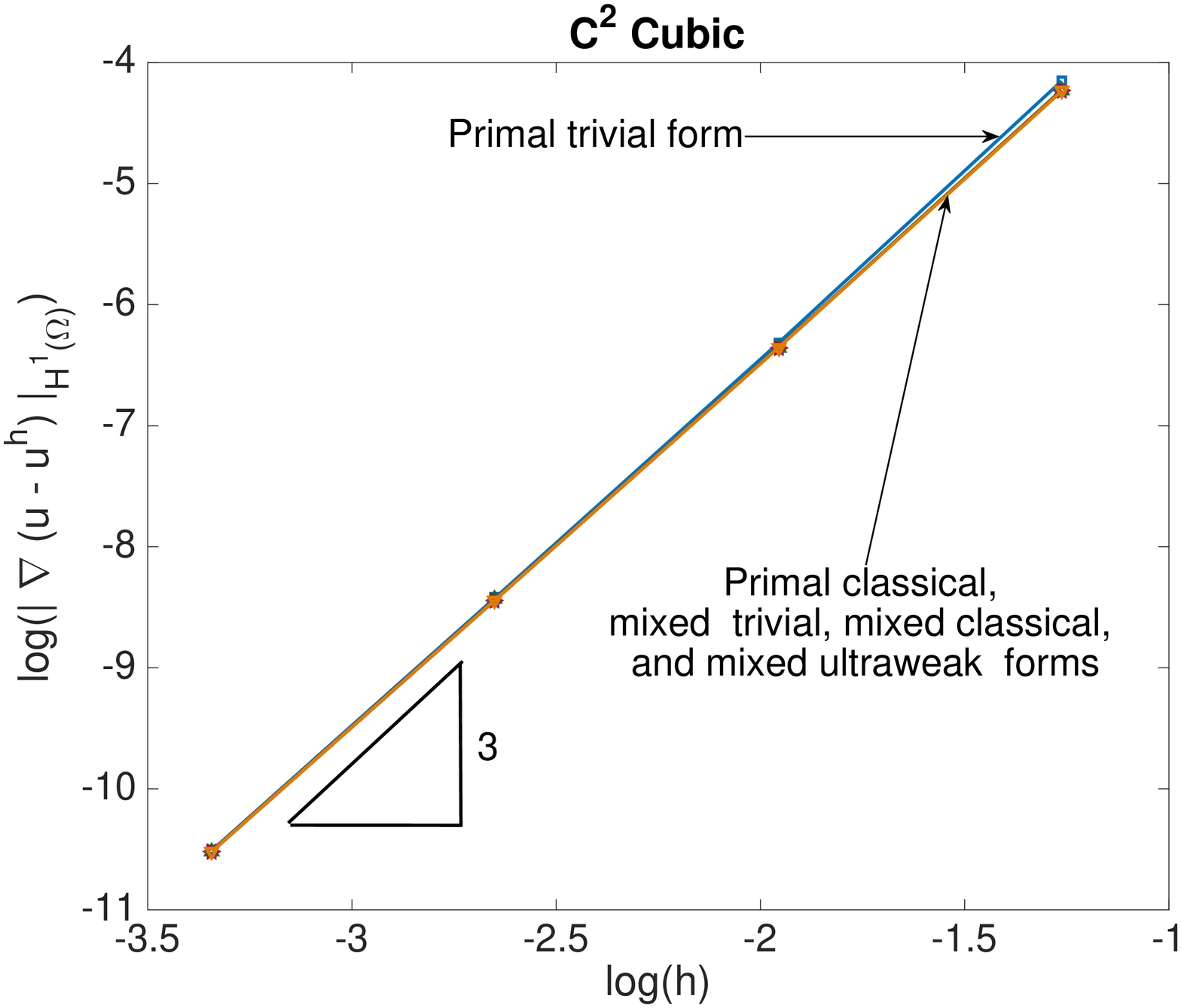}
\centering\includegraphics[width=6.0cm]{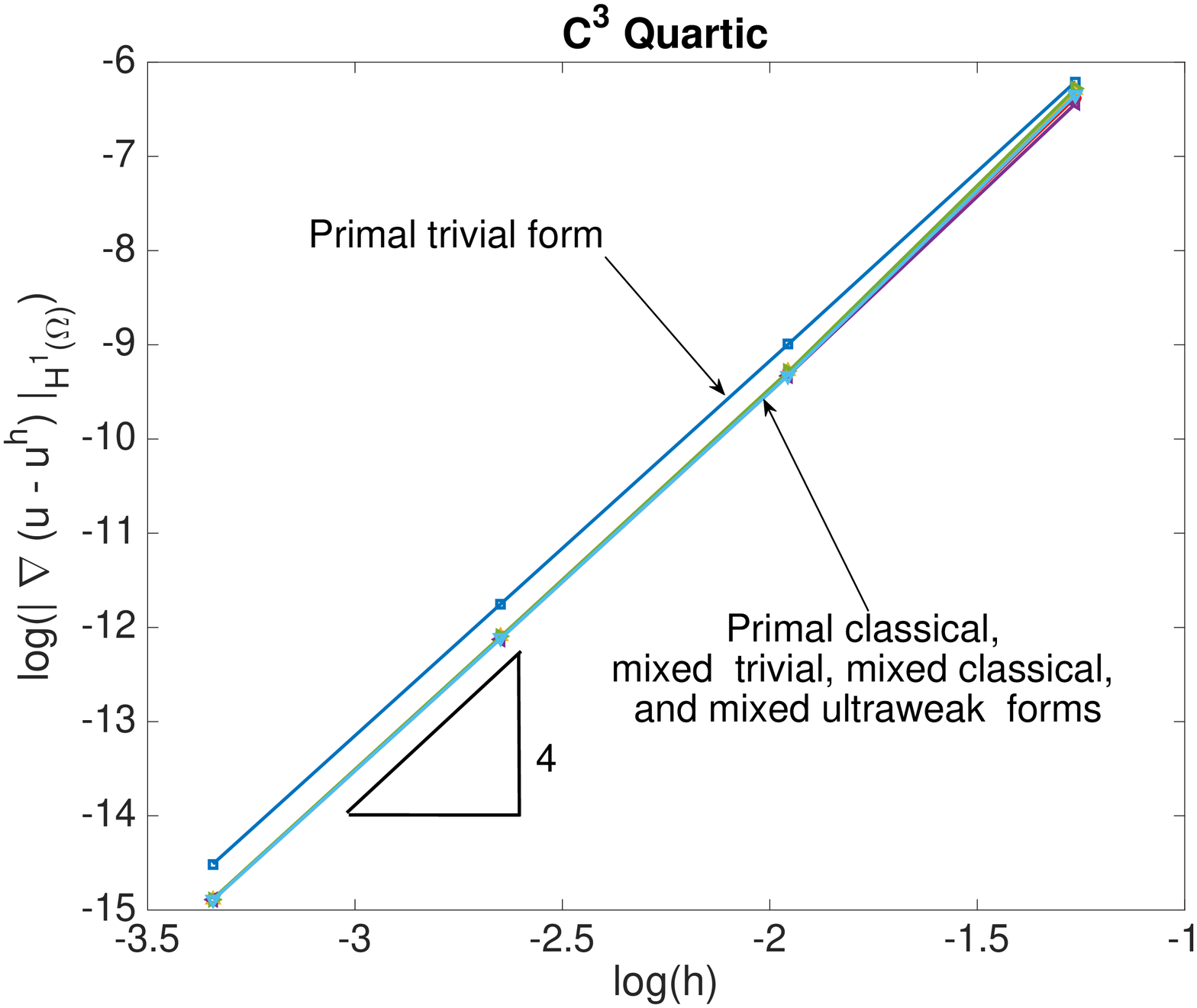}
\centering\includegraphics[width=6.0cm]{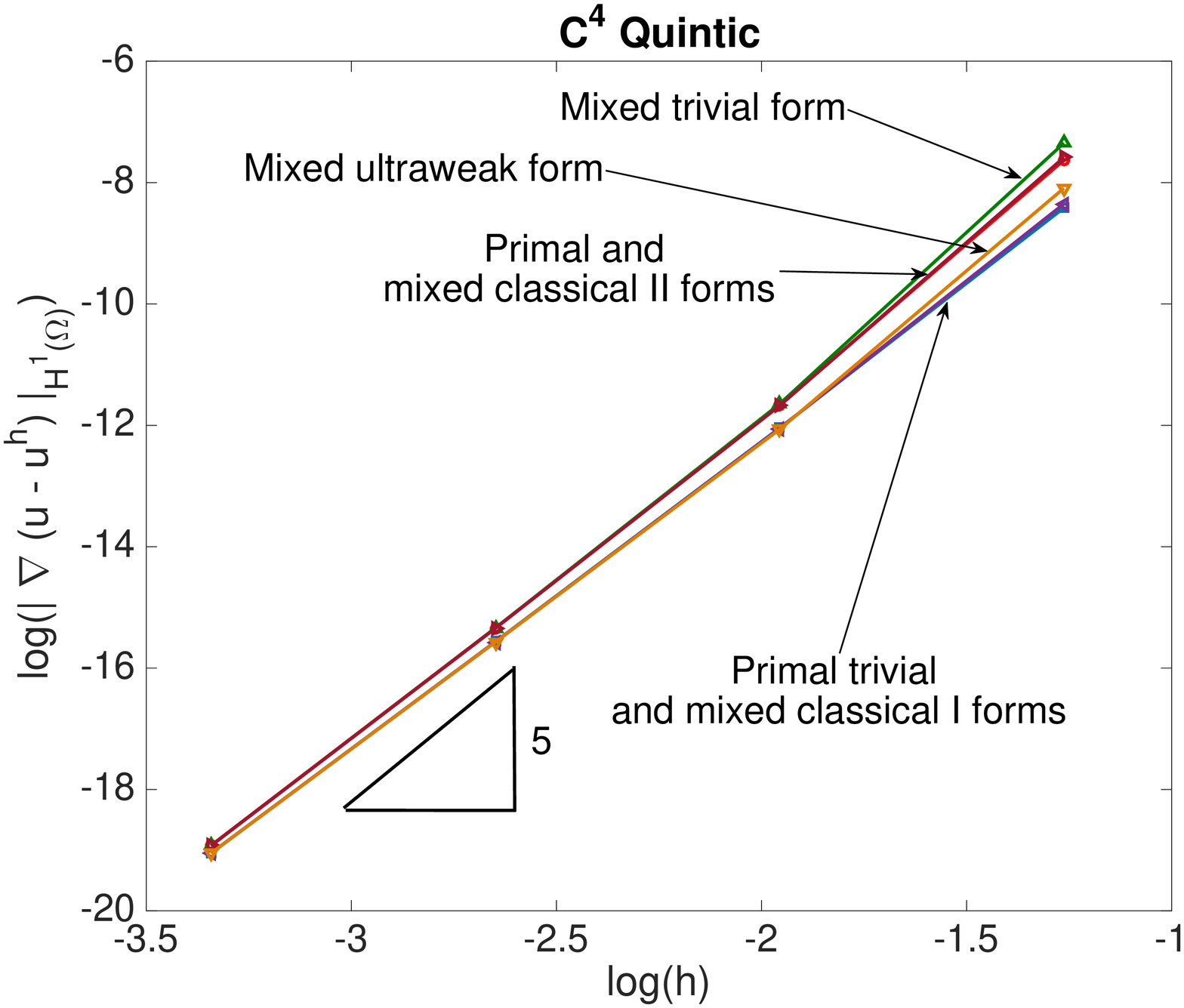}
\caption{$H^1$ semi-norm errors when using six formulations with isogeometric elements.}
\label{fig:p2345uh1pm}
\end{figure}

Both the primal and mixed trivial formulations do not involve any integration by parts in their variational formulations and their test functions are in $L^2$. The difference is that the primal formulation results from \eqref{eq:pde} while the mixed formulation results from \eqref{eq:pdee}.  We use $C^{-1}$ basis functions for the test spaces for these formulations. Consequently, we obtain a matrix $G$ in \eqref{eq:mp} which is a block-diagonal matrix. The independence of each elemental block from the rest allows these methods to be computable efficiently (cheap elemental inversion) and thus relevant for practical purposes. All other formulations involve integration by parts to pass derivatives to the test functions, which in return results in a matrix $G$ in \eqref{eq:mp} which is not block-diagonal. This makes these formulations interesting from the theoretical point of view, but untractable in most general meshes, even though splitting schemes make these methods viable on tensor-product meshes \cite{los2018}.  Thus, once we show that all these formulations deliver equivalent results at the discrete level, we focus on these strong/trivial formulations as they are computationally advantageous. We chose not to compare against well established DPG technologies to brake test spaces, as the goal is simply to show the equivalence of the different variational forms rather than derive alternative computable methods. Lastly, to compare the primal trivial formulation with the first-order system least-square (FOSLS) method \cite{cai1994first}, the difference is that the primal trivial formulation does not lead to first order system and it introduces a Gramm matrix to solve for the residual errors simultaneously.

\begin{figure}[ht]
\centering\includegraphics[width=6.0cm]{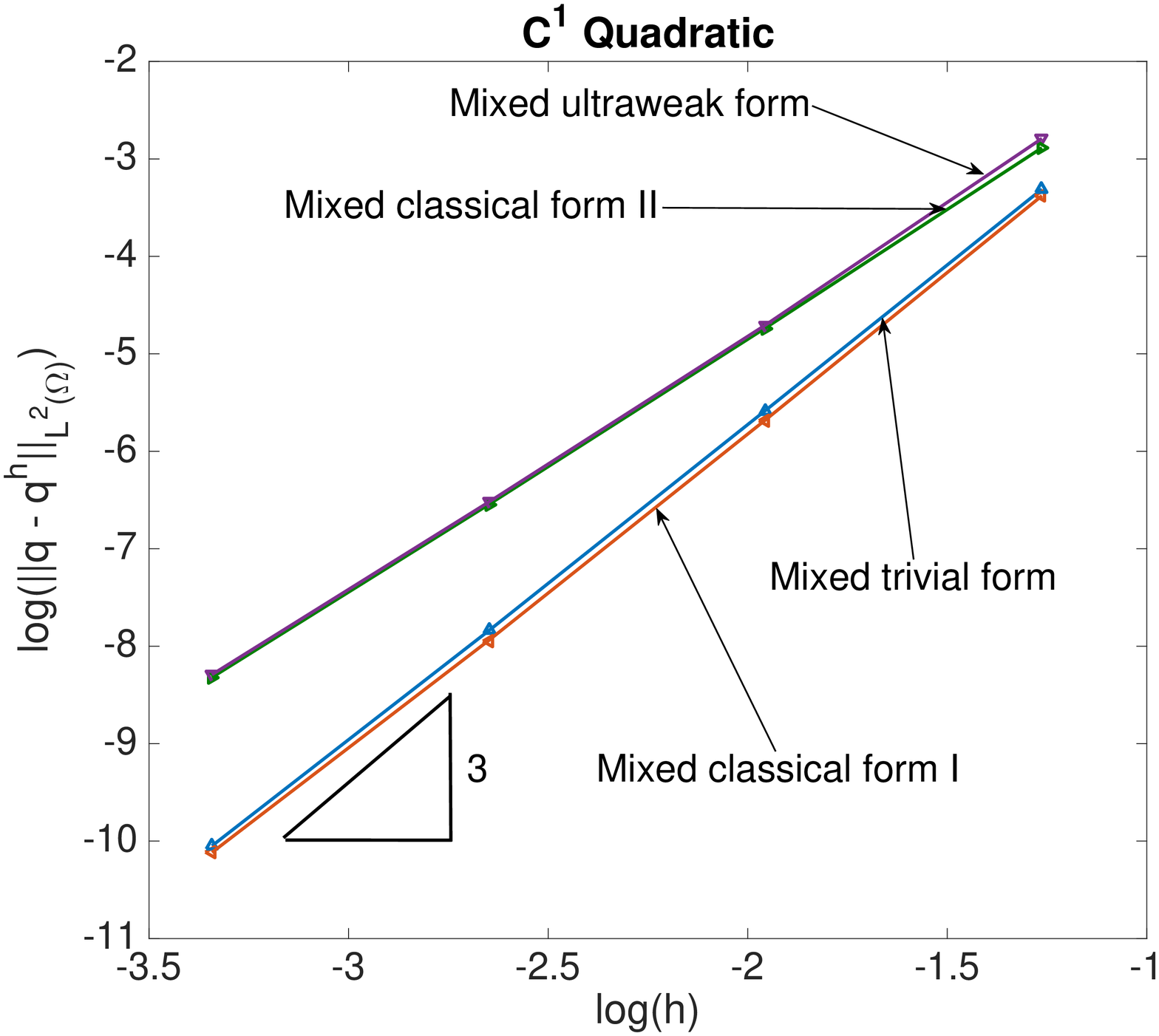} 
\centering\includegraphics[width=6.0cm]{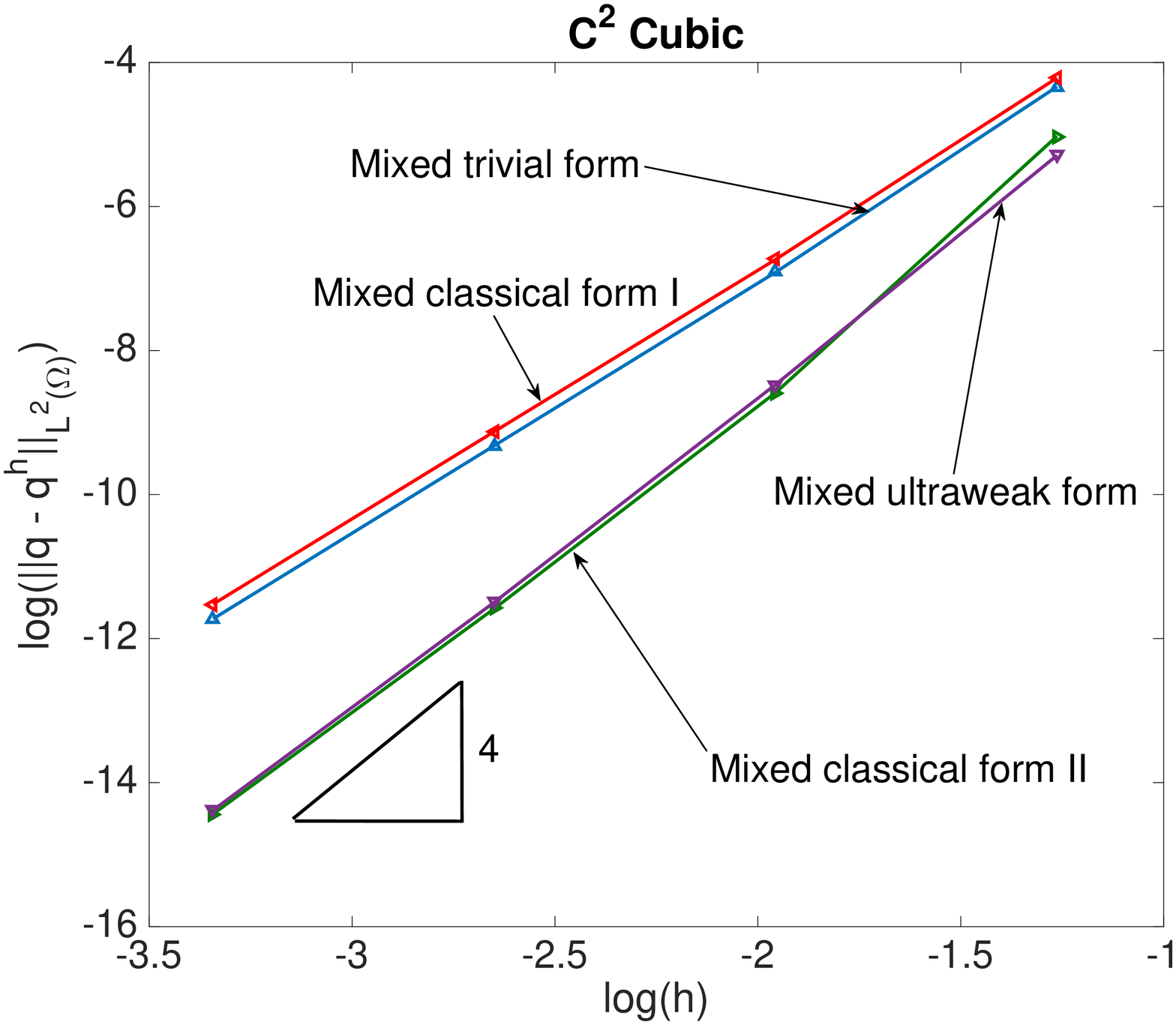}
\centering\includegraphics[width=6.0cm]{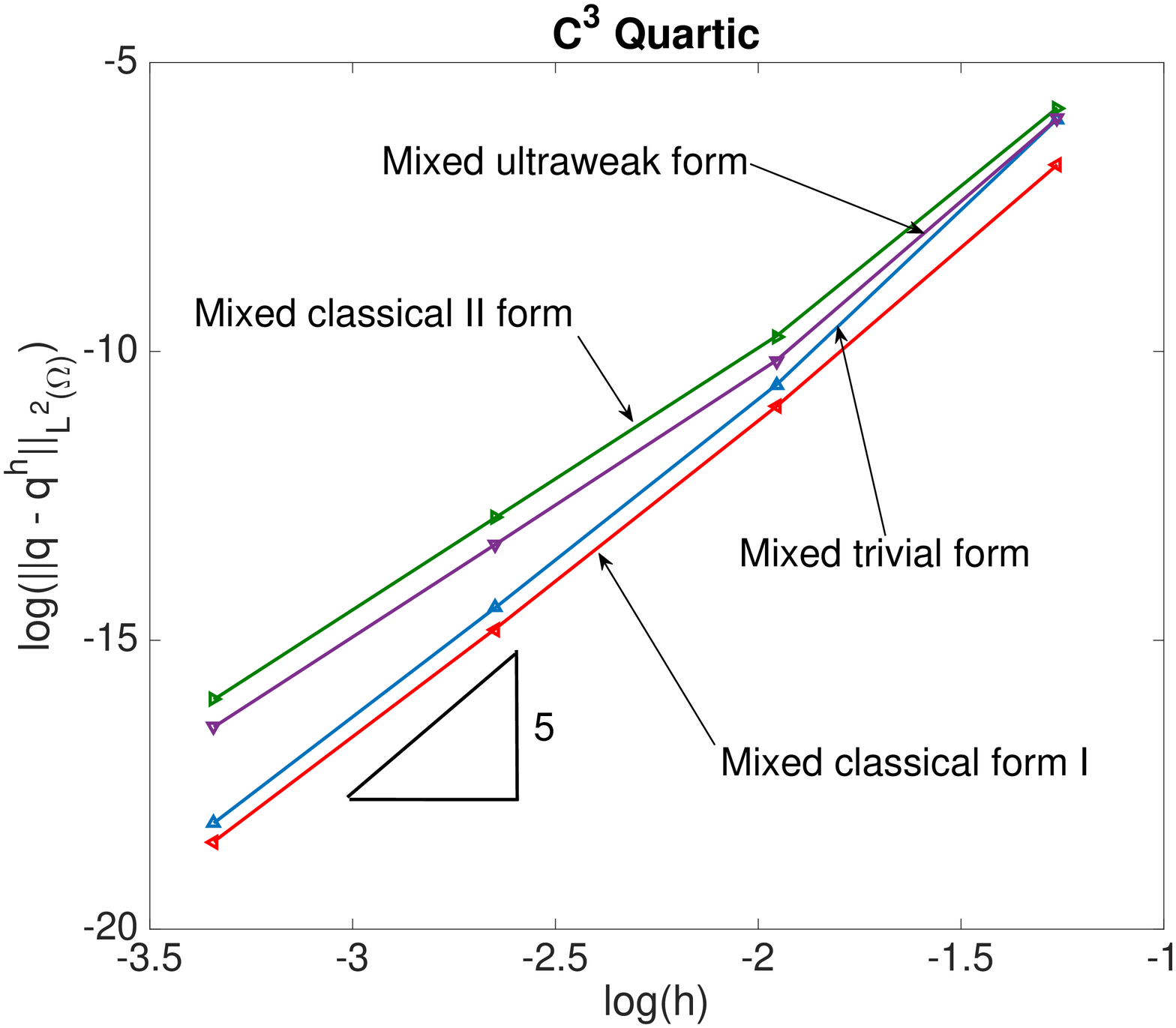}
\centering\includegraphics[width=6.0cm]{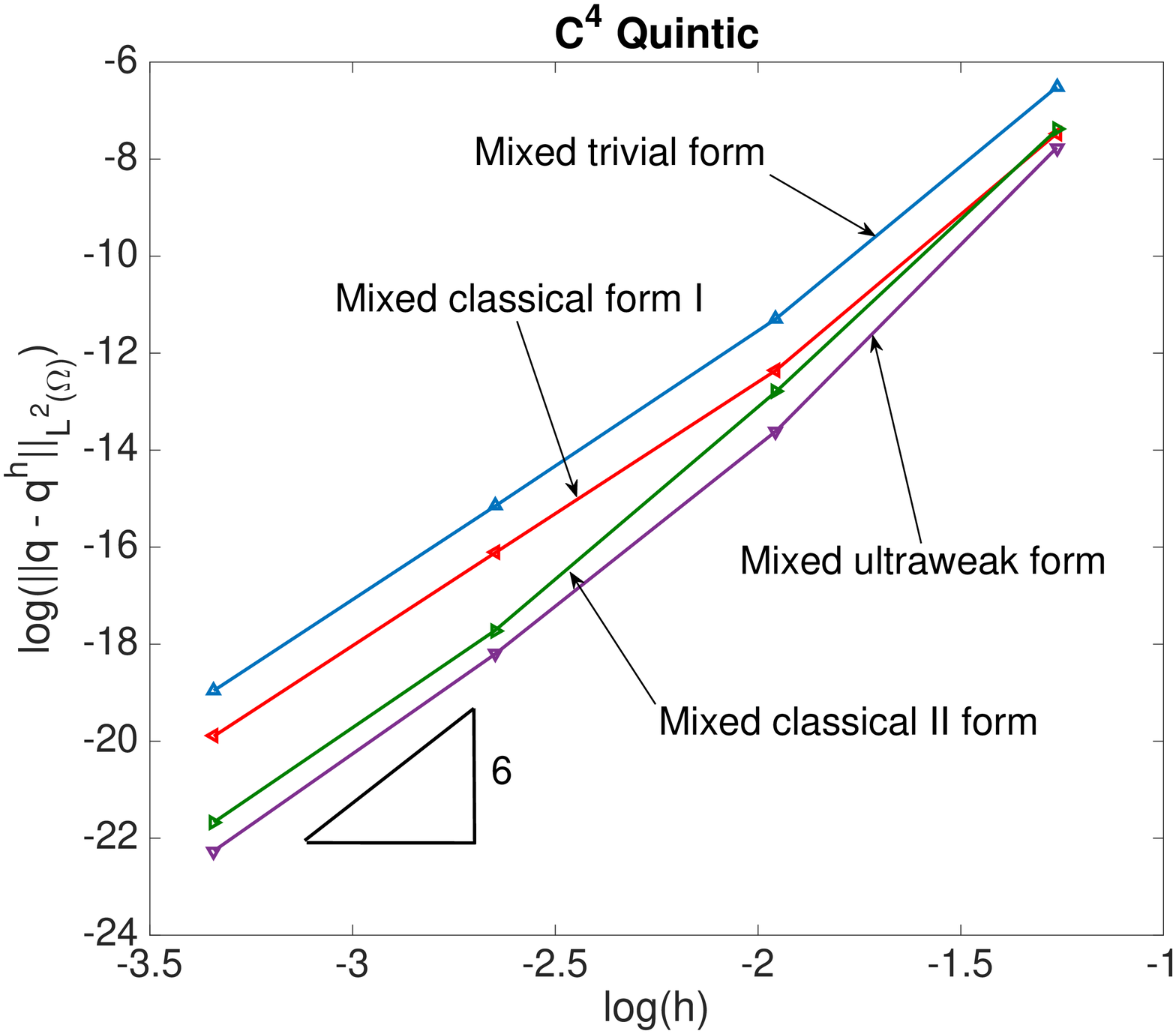}
\caption{$L^2$ errors of $\bfs{q}^h$ when using mixed formulations with isogeometric elements.}
\label{fig:p2345ql2m}
\end{figure}

Figure \ref{fig:p2345uh1pm} shows the $H^1$ semi-norm errors when using all six formulations for $p=2,3,4,5$. For the mixed formulations, we approximate the flux $\bfs{q}^h$ using basis functions of the same polynomial order as for the solution $u^h$. The parameters of the bilinear form $g(\cdot, \cdot)$ are set to be the default values.  The mesh configurations are $5\times5, 10\times10, 20\times20, 40\times40$. Herein, we plot the errors in natural logarithmic scale. As predicted, in all scenarios, the $H^1$ semi-norm errors converge in optimal rates, that is, order $p$ for all the variational formulations. This confirms the theoretical result \eqref{eq:uh1err}. Therefore, all formulations are equivalent in the quality of the approximation they deliver. Interestingly, all primal trivial forms for even orders have optimal convergence rates, but the constants seem to be worse for this discretization for even orders than all the other weak forms we compared against. Nevertheless, for odd order polynomials both the rate of convergence and constant are comparable to those observed for all the other variational forms.

\begin{figure}[!ht] 
\centering\includegraphics[width=6.0cm]{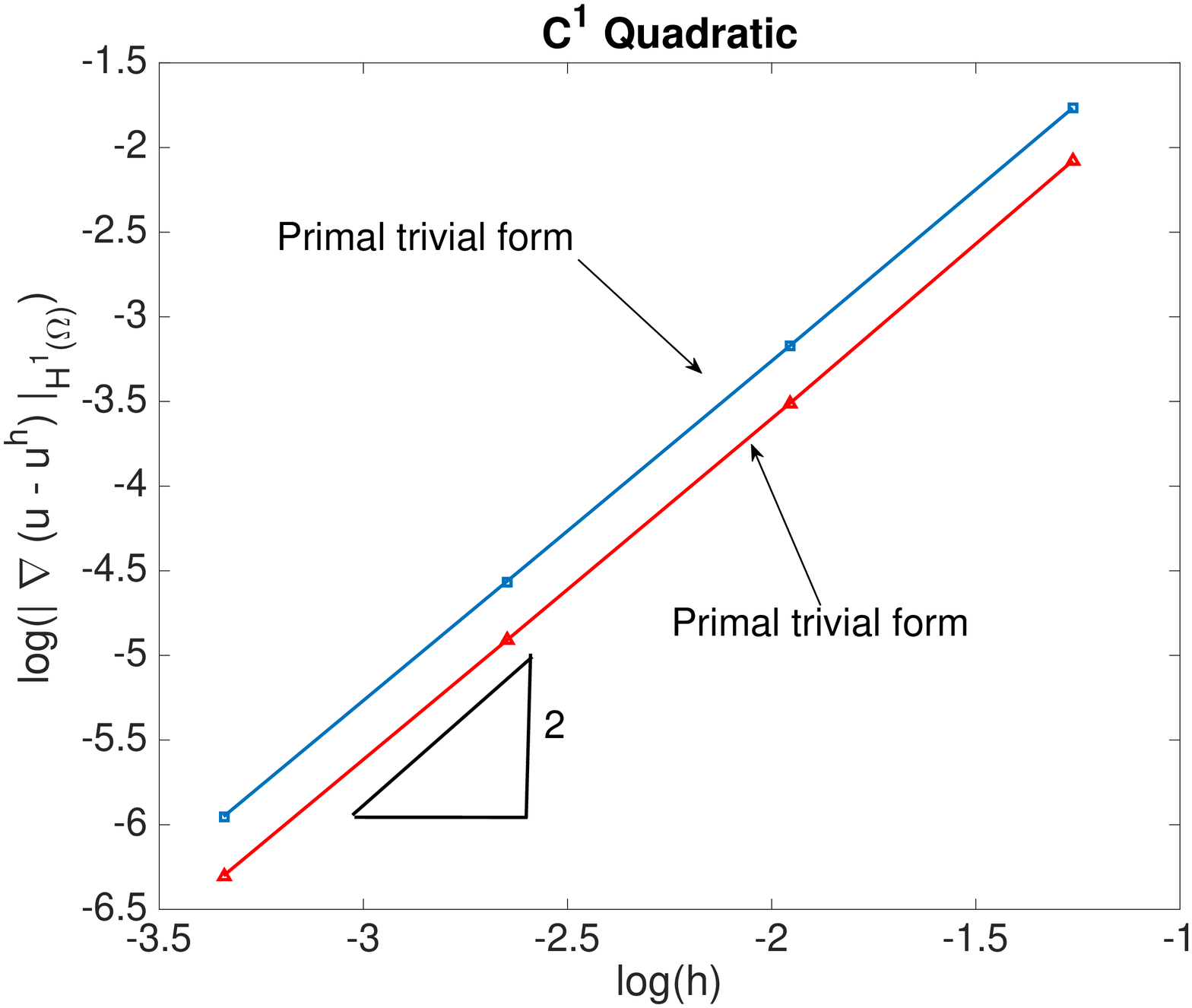} 
\centering\includegraphics[width=6.0cm]{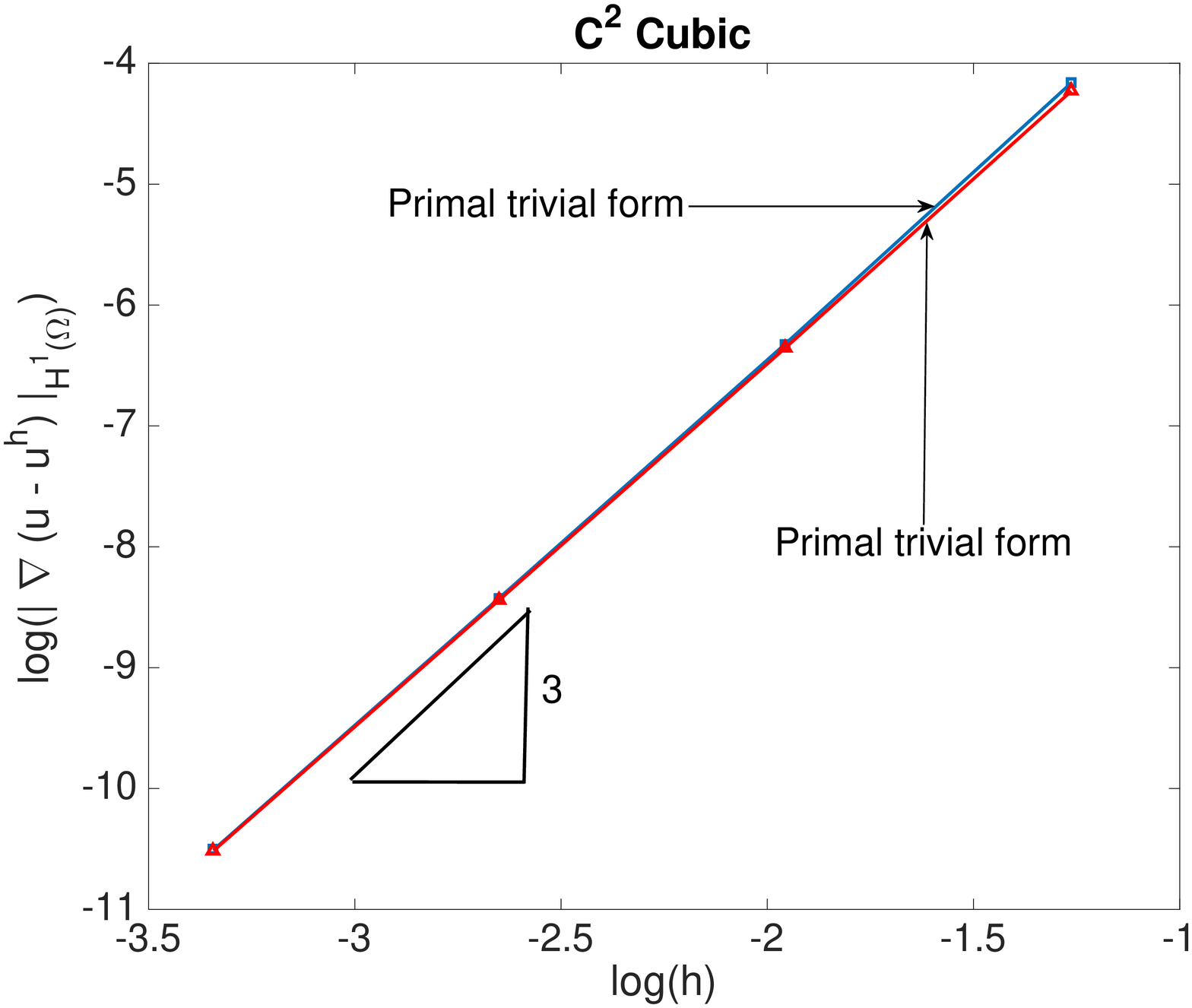}
\caption{$H^1$ semi-norm errors of $u^h$ when using primal trivial formulation with isogeometric elements and an $L^2$ product defining the Gramm matrix.}
\label{fig:p23auh1}
\end{figure}

The mixed formulations introduce auxiliary variables for the fluxes. Thus, for the same mesh configuration, the resulting matrix systems of the mixed formulations have dimensions which are three times larger than the primal ones. However, the mixed formulations deliver more accurate fluxes approximations.

Figure \ref{fig:p2345ql2m} shows the flux errors in the $L^2$ norm when using the mixed forms for $p=2,3,4,5.$ The fluxes are of optimal convergence rates $p+1$, which confirms the error estimate \eqref{eq:qhl2}. For $p=2,4$, we observe super-convergence rates approximately $p+2$. Once again, these results  numerically verify that these four formulations are equivalent. Therefore, for the following numerical results, we only show the results for the strong primal and mixed formulations.

We also study the formulations with different choices of the auxiliary bilinear forms $g(\cdot, \cdot)$ in \eqref{eq:ra} and \eqref{eq:ma}.  Figure \ref{fig:p23auh1} shows the numerical errors $u^h$ in $H^1$ semi-norm when using the inner product $g(v, w) = (v, w)_\Omega$ in \eqref{eq:ra} for the trivial primal form and $g((v, \bfs{r}), (w, \bfs{p}) )  = (v, w)_\Omega + (\bfs{r}, \bfs{p})_\Omega$ in \eqref{eq:ma} for the trival mixed form. Other inner products are also possible, which is in agreement with the DPG methodology.  The convergence rates in the $H^1$ semi-norm of all these scenarios are optimal.

\section{Concluding remarks} \label{sec:con}
We introduce a residual minimization based mixed formulation for solving partial differential equations. A key feature of this method is the framework which uses highly-continuous B-splines for the trial spaces and basis functions with minimal regularity and lower order polynomials for the test spaces. The method shares with the discontinuous Petrov-Galerkin methodology the idea of stabilizing the formulation considering an adequate norm for the test space and unifies the interpretation of several methods such as the classical finite element method, isogeometric analysis, and isogeometric collocation methods. Under the standard assumption, the proposed variational formulations are stable and result in optimal approximation properties.

\section*{Acknowledgement}
This publication was made possible in part by the CSIRO Professorial Chair in Computational Geoscience at Curtin University and the Deep Earth Imaging Enterprise Future Science Platforms of the Commonwealth Scientific Industrial Research Organisation, CSIRO, of Australia. The J. Tinsley Oden Faculty Fellowship Research Program at the Institute for Computational Engineering and Sciences (ICES) of the University of Texas at Austin has partially supported the visits of VMC to ICES. 

%
%
%
\bibliographystyle{splncs04}
%

\end{document}